\documentclass{amsart}
\usepackage[english]{babel}

\usepackage{amsmath}
\usepackage{amsthm}
\usepackage{amsfonts}
\usepackage{amssymb}
\usepackage{euscript}
\usepackage[all]{xy}
\newcommand{\dimv}{\underline{\dim}}

\newtheorem{thm}{Theorem}[section]

\newtheorem{lem}[thm]{Lemma}
\newtheorem{cor}[thm]{Corollary}

\newtheorem{question}[thm]{Question}

\begin{document}
\title{Embeddings of representations}
\author{Kathrin Kerkmann, Markus Reineke}
\address{Kathrin Kerkmann:\newline
Fakult\"at f\"ur Mathematik, Ruhr-Universit\"at Bochum, D - 44780 Bochum}
\email{kathrin.kerkmann@ruhr-uni-bochum.de}
\address{Markus Reineke:\newline
Fachbereich C - Mathematik, Bergische Universit\"at Wuppertal, D - 42097 Wuppertal}
\email{mreineke@uni-wuppertal.de}
\begin{abstract} We derive ``numerical'' criteria for the existence of embeddings of representations of finite dimensional algebras.\end{abstract}

\date{}
\maketitle
\parindent0pt
\section{Introduction}
By a classical result of M. Auslander \cite{Aus,BoAus}, a finite dimensional representation $M$ of a finitely generated algebra $A$ is determined up to isomorphism by the dimensions of homomorphism spaces to it, that is, two such representations $M$ and $N$ of $A$ are isomorphic if and only if $\dim{\rm Hom}(U,M)=\dim{\rm Hom}(U,N)$ for all (indecomposable) representations $U$ of $A$.\\[1ex]
In light of this fact, one can ask for ``numerical'' criteria for representation-theoretic properties. One example is the characterization of degenerations \cite{Bo} $M\leq_{\deg}N$ of representations of algebras of finite representation type by the condition $\dim{\rm Hom}(U,M)\leq\dim{\rm Hom}(U,N)$ for all $U$ \cite{Zw}.\\[1ex]
The aim of the present paper is to prove numerical criteria for situations related to embeddings of representations. This question is motivated by a study of quiver Grassmannians for representations of Dynkin quivers, for which specific geometric properties can be expected (in contrast to arbitrary quiver Grassmannians, see \cite{RQG}). The first step in this direction is a criterion for nonemptyness of a quiver Grassmannian, which will be proven in the following form:

\begin{thm}\label{gr1} A representation $M$ of a Dynkin quiver $Q$ with associated Euler form $\langle\_,\_\rangle$ admits a subrepresentation of dimension vector ${\bf e}$ if and only if $\dim{\rm Hom}(U,M)\geq\langle{\bf dim}\ U,{\bf e}\rangle$ for all (indecomposable) representations $U$ of $Q$.
\end{thm}

Not directly related, but in the same spirit, we find a quite general sufficient criterion for irreducibility of a Dynkin quiver Grassmannian:
\begin{thm}\label{gr2} Given a dimension vector ${\bf e}$ and a representation $M$ as before such that the following inequalities hold:
\begin{enumerate}
\item $\dim{\rm Hom}(M,U)\leq\langle{\bf e},{\bf dim}\, U\rangle$ for all non-injective indecomposable $U$,
\item $\dim{\rm Hom}(U,M)\leq\langle{\bf dim}\, U,{\bf dim}\, M-{\bf e}\rangle$ for all non-projective indecomposables $U$,
\end{enumerate}
the quiver Grassmannian ${\rm Gr}_{\bf e}(M)$ is irreducible of dimension $$\dim{\rm Gr}_{\bf e}(M)=\langle{\bf e},{\bf dim}\, M-{\bf e}\rangle.$$
\end{thm}

 Both results were predicted by extensive numerical experiments for a type $A_3$ quiver in the first named author's master thesis \cite{K}.\\[1ex]
The other main topic of this paper concerns the much finer problem to numerically characterize embeddings between two given representations. In this direction, we prove




\begin{thm}\label{thm1} Let $A$ be an arbitrary finite dimensional algebra over an algebraically closed field $k$, and let $M$ and $N$ be finite dimensional representations of $A$. Then the following are equivalent:
\begin{enumerate}
\item For all large enough $r\geq 1$, there exists an embedding $N^r\rightarrow M^r$,
\item for all surjections $U\rightarrow V$ of representations of $A$, we have
\begin{equation}\label{nc2}[U,N]-[V,N]\leq [U,M]-[V,M],\end{equation}
\item the estimate (\ref{nc2}) holds for all quotients $U=N^k\rightarrow N^k/S=V$, where $S$ is a simple subrepresentation of $N^k$ and $k\leq\dim{\rm Hom}(S,N)$.
\end{enumerate}
\end{thm}

Note that the numerical condition (\ref{nc2}) is insensitive to multiplicities, so that one cannot expect to characterize existence of an actual embedding $N\subset M$ in general. In fact, in Section \ref{section_ex}, we will exhibit a (low-dimensional) example of representations $N$ and $M$ over the three-arrow Kronecker quiver such that $N^2$ can be embedded into $M^2$, but $N$ cannot be embedded into $M$. It is rather natural to ask for which algebras the conditon (\ref{nc2}) already characterizes embeddings $N\subset M$; at least this holds for an equioriented type $A$ quiver, see Section \ref{section_ex}.\\[1ex]
Theorem \ref{thm1} will be proved in Section \ref{s2} using a linear algebra lemma which will be derived in Section \ref{app} from a slight generalization of a theorem of W. Crawley-Boevey. We complement this result by examples and remarks in Section \ref{section_ex}. Theorems \ref{gr1} and \ref{gr2} will be proved in Section \ref{s4} and Section \ref{s5}, respectively.

\section{Proof of Theorem \ref{thm1}}\label{s2}

In this section we prove Theorem \ref{thm1} stated above using a linear algebra result, Lemma \ref{lal}, which will be proved in Section \ref{app}.\\[1ex]
The third statement of Theorem \ref{thm1} is a special case of the second one. The first statement implies the second as follows:\\[1ex]
for the exact sequences
$$0\rightarrow N^r\rightarrow M^r\rightarrow C\rightarrow 0$$
and
$$0\rightarrow K\rightarrow U\rightarrow V\rightarrow 0$$
we construct a diagram of spaces of homomorphisms (abbreviating ${\rm Hom}(X,Y)$ by $(X,Y)$):
$$\begin{array}{ccccccc}
&&0&&0&&0\\
&&\downarrow&&\downarrow&&\downarrow\\
0&\rightarrow&(V,N^r)&\rightarrow&(V,M^r)&\rightarrow&(V,C)\\
&&\downarrow&&\downarrow&&\downarrow\\
0&\rightarrow&(U,N^r)&\rightarrow&(U,M^r)&\rightarrow&(U,C)\\
&&\downarrow&&\downarrow&&\downarrow\\
0&\rightarrow&(K,N^r)&\rightarrow&(K,M^r)&\rightarrow&(K,C),
\end{array}$$
resulting in a left exact sequence
$$0\rightarrow(V,N^r)\rightarrow(V,M^r)\oplus(U,N^r)\rightarrow(U,M^r).$$
The estimate (\ref{nc2}) follows. To prove that the third statement implies the first, we borrow the following result from the appendix:

\begin{lem}\label{lal} Let $Z\subset{\rm Hom}(V,W)$ be a subspace of linear maps between finite dimensional $k$-vector spaces $V$ and $W$. Assume that $\dim Z(U)\geq\dim U$ for all $U\subset V$. Then, for all large enough $r$, there exists a matrix $F\in M_{r\times r}(Z)$ defining an injective map from $V^r$ to $W^r$.
\end{lem}

Assume given $M$ and $N$ as in the theorem such that (iii) holds. We will show the following: for every isotypical component $I$ of the socle of $N$, there exists, for large enough $r_I$, a map $f_I$ from $N^{r_I}$ to $M^{r_I}$ whose kernel has zero intersection with $I^{r_I}$. Since there are only finitely many isotypical components, this implies that, for $r$ large enough, a generic map from $N^r$ to $M^r$ is non-vanishing on the socle of $N^r$ and thus injective, proving the theorem.\\[1ex]
So assume $I$ is an isotypical component of $N$, thus $I\simeq S^{[S,N]}$ for a simple $S$, where $[S,N]=\dim{\rm Hom}(S,N)$. We would like to apply Lemma \ref{lal} and define $V={\rm Hom}(S,N)$, $W={\rm Hom}(S,M)$ and $Z$ as the image of the natural map
$${\rm Hom}_A(N,M)\rightarrow{\rm Hom}_k({\rm Hom}_A(S,N),{\rm Hom}_A(S,M))$$
given by composition. We show that the assumption of Lemma \ref{lal} is fulfilled: let $U$ be a subspace of ${\rm Hom}(S,N)$, spanned by, say, $g_1,\ldots,g_l$. We consider the exact sequence
$$0\rightarrow S\stackrel{G}{\rightarrow} N^l\rightarrow\bar{N}\rightarrow 0,$$
where $G=[g_1,\ldots,g_l]$. Applying ${\rm Hom}(\_,N)$ we get an exact sequence
$$0\rightarrow{\rm Hom}(\bar{N},N)\rightarrow{\rm Hom}(N^l,N)\stackrel{\alpha}{\rightarrow}{\rm Hom}(S,N).$$
By definition, the image of $\alpha$ contains $U$. On the other hand, we have an exact sequence
$$0\rightarrow{\rm Hom}(\bar{N},M)\rightarrow{\rm Hom}(N^l,M)\stackrel{\gamma}{\rightarrow}{\rm Hom}(S,M),$$
and the image of $\gamma$ by definition equals $\sum_{i=1}^l{\rm Hom}(N,M)g_i$, which in the above notation can be rewritten as $Z(U)$.\\[1ex]
By the assumption (iii) applied to the surjection $N^l\rightarrow\bar{N}$, we have the following estimate:
$$\dim U\leq\dim{\rm Im}(\alpha)=[N^l,N]-[\bar{N},N]\leq$$
$$\leq[N^l,M]-[\bar{N},M]=\dim{\rm Im}(\gamma)=\dim Z(U).$$
Applying Lemma \ref{lal}, we conclude existence of a map $F:N^r\rightarrow M^r$ for large enough $r\geq 1$ such that the induced map (given by composition) from ${\rm Hom}(S,N^r)$ to ${\rm Hom}(S,M^r)$ is injective. This implies that the intersection of the kernel of $F$ with $I^r$ is zero.

\section{Examples and remarks}\label{section_ex}

The guiding question for this section is:

\begin{question}\label{question} For which algebras $A$ does the condition (\ref{nc2}) of Theorem \ref{thm1} already characterize existence of an embedding of $N$ into $M$?
\end{question}

Let us first verify that a path algebra of an equioriented type $A$ quiver $1\rightarrow\ldots\rightarrow n$ belongs to this class of algebras. Denote by $U_{i,j}$ for $1\leq i\leq j\leq n$ the indecomposable supported on the interval $[i,j]$. The conditions (\ref{nc2}) for two given representations $N$ and $M$, applied to the surjections $U_{i,j+1}\rightarrow U_{i,j}$, read $$\sum_{k\leq i}n_{k,j}\leq\sum_{k\leq i}m_{k,j}$$ for all $i\leq j$, where $n_{i,j}$ (resp.~$m_{i,j}$) denotes the multiplicity of $U_{i,j}$ as a direct summand of $N$ (resp.~$M$). Given these conditions, an embedding $N\rightarrow M$ can be constructed already as an appropriate direct sum of the embeddings $U_{i+1,j}\rightarrow U_{i,j}$.\\[1ex]
Now we turn to counterexamples. Let $Q$ be the three-arrow Kronecker quiver, that is, $Q$ has vertices $i$ and $j$ and three arrows from $i$ to $j$. Let $P_i$ be the indecomposable projective attached to the source $i$, which is of dimension vector $(1,3)$. Let $M$ be the representation of dimension vector $(3,3)$ such that the arrows are represented by the following matrices:
$$\left(\begin{array}{rrr}1&0&0\\ 0&0&0\\ 0&0&-1\end{array}\right),\; \left(\begin{array}{rrr}0&0&0\\ 1&0&0\\ 0&1&0\end{array}\right),\;\left(\begin{array}{rrr}0&1&0\\ 0&0&1\\ 0&0&0\end{array}\right)$$
(note that it has a quiver of type $\tilde{A}_5$ as coefficient quiver).\\[1ex]
Then $P_i$ does not embed into $M$: in a map $f$ from $P_i$ to $M$, if $f_i$ is given by a vector with coordinates $a$, $b$ and $c$, the map $f_j$ is necessarily given by the matrix
$$\left(\begin{array}{rrr} a&0&b\\ 0&a&c\\ -c&b&0\end{array}\right),$$
which has zero determinant. But $P_i^2$ embeds into $M^2$, for example via the map $g$ given by the matrices
$$g_i=\left(\begin{array}{rr}0&1\\ 1&0\\ 0&0\\ 1&0\\ 0&0\\ 0&1\end{array}\right),\; g_j=\left(\begin{array}{rrrrrr}0&0&1&1&0&0\\ 0&0&0&0&1&0\\ 0&1&0&0&0&0\\ 1&0&0&0&0&0\\ 0&1&0&0&0&1\\ 0&0&0&-1&0&0\end{array}\right)$$
 (note that $g_j$ has determinant $1$).\\[1ex]
Via the translation to a problem about (semistable) representations of generalized Kronecker quivers in the appendix, this example was found using the description of semiinvariants for representations of $Q$ of dimension vector $(3,3)$ in \cite{Do}; namely, the semiinvariant
$$(A,B,C)\mapsto\det\left(\begin{array}{rr}B&A\\ A&C\end{array}\right)$$
of triples of $3\times 3$-matrices is algebraically independent from polarizations of the determinant, that is, from the semiinvariants arising as coefficients of $x$, $y$ and $z$ in $\det(xA+yB+zC)$.\\[1ex]
Next we remark that assumptions on the ground field play an essential role for Question \ref{question}. Consider the path algebra of the following quiver of type $D_4$:
$$\begin{array}{ccccc}2&\leftarrow&1&\rightarrow&4\\ &&\downarrow&&\\ &&3&&\end{array}$$
Then, for any field $k$ with at least three elements, the indecomposable projective $P_1$ embeds into the five-dimensional indecomposable $X$, which is not true for $k={\bf F}_2$. But the numerical condition (\ref{nc2}) is always fulfilled due to independence of the representation theory of a Dynkin quiver from the ground field. As remarked by C.~M.~Ringel, it is also easy to see that, similarly to the previous example, $P_1^2$ admits an embedding into $X^2$ even over ${\bf F}_2$.\\[1ex]
Passing from $A$ to $A^{\rm op}$, we have an obvious dual statement to Theorem \ref{thm1} characterizing the existence of surjections. Namely, for given representations $U$ and $V$, there exists a surjection $U^r\rightarrow V^r$ for large enough $r\geq 1$ if and only if
$$[U,N]-[V,N]\leq [U,M]-[V,M]$$
for all embeddings $N\subset M$.\\[1ex]
Now we reinterprete Theorem \ref{thm1} and Question \ref{question} in the language of cones in split Grothendieck groups. Let $\Lambda$ be a lattice, and let $C\subset\Lambda$ be a cone. $C$ is called saturated if $rx\in C$ for some $r\geq 1$ inplies $x\in C$. The saturation $\hat{C}$ of $C$ is defined as the set of all $x\in\Lambda$ such that $rx\in C$ for some $r\geq 1$. Define ${\rm Inj}$ as the set of all pairs $(N,M)$ in the product $K_0(\bmod A)\times K_0(\bmod A)$ of split Grothendieck groups such that $N$ embeds into $M$. Then Question \ref{question} is equivalent to asking for which algebras $A$ this cone is saturated.\\[1ex]
A somewhat related (but much deeper) saturation result is the following: let $A=k[[t]]$ be the ring of formal power series; thus isomorphism classes of finite dimensional representations of $A$ are parametrized by partitions. Under this identification, the cone of triples $(N,X,N)$ in $K_0(\bmod\, A)^3$ fitting into a short exact sequence $0\rightarrow N\rightarrow X\rightarrow N\rightarrow 0$ equals the Littlewood-Richardson cone by \cite{Mac}, which is proved to be saturated in \cite{KT}.\\[1ex]
Given lattices $\Lambda$ and $\Lambda'$ and a (biadditive) pairing $(\_,\_):\Lambda\times\Lambda'\rightarrow{\bf Z}$, we call cones $C\subset\Lambda$ and $C'\subset\Lambda'$ dual if 
\begin{enumerate}
\item $x\in C$ iff $(x,C')\geq 0$,
\item $y\in C'$ iff $(C,y)\geq 0$.
\end{enumerate}
Define ${\rm Surj}$ as the cone of all pairs $(U,V)\in K_0(\bmod A)\times K_0(\bmod A)$ admitting a surjection $U\rightarrow V$. We define a pairing on 
$K_0(\bmod A)\times K_0(\bmod A)$ by
$$((N,M),(U,V))\mapsto [U-V,M-N]$$
(note that the closely related pairing $(X,Y)\mapsto[X,Y]$ on $K_0(\bmod A)$ is nondegenerate by Auslander's theorem cited in the introduction). Using this notation, we can reinterprete Theorem \ref{thm1} (and its dual version stated above) as follows:
\begin{cor} The cones $\widehat{\rm Inj}$ and $\widehat{\rm Surj}$ are dual.
\end{cor}

\section{Non-emptyness of quiver Grassmannians}\label{s4}

In this section, we assume $Q$ to be a Dynkin quiver with set of vertices $Q_0$ and Euler form $\langle\_,\_\rangle$ and recollect some known results on generic properties of representations of $Q$.\\[1ex]
For a dimension vector ${\bf e}\in{\bf N}Q_0$, let $G_{\bf e}$ be the unique (up to isomorphism) representation of dimension vector ${\bf e}$ with vanishing group of selfextensions, or, in other words, the unique representation with an open orbit in its variety of representations. Using this notation, we can reformulate a theorem of A. Schofield \cite{S} as: $G_{\bf d}$ admits a subrepresentation of dimension vector ${\bf e}$ if and only if ${\rm Ext}^1(G_{\bf e},G_{{\bf d}-{\bf e}})=0$.\\[1ex]
We also recall the concept of generic extensions from \cite{R}: given representations $M$ and $N$ of $Q$, there exists a unique (up to isomorphism) representation $M*N$ which is an extension of $M$ by $N$, and has minimal dimension of its endomorphism ring with this property; in other words, it has a dense orbit in the subvariety (of the corresponding representation variety) consisting of all extensions of $M$ by $N$. The representation $M*N$ is called the generic extension of $M$ by $N$.\\[1ex]
The following is proven in \cite{R}: if $M$ degenerates to $M'$ and $N$ degenerates to $N'$, and $X$ is an extension of $M'$ by $N'$, then $M*N$ degenerates to $X$. Using this result and the fact that there is no proper degeneration to a representation $G_{\bf e}$, we see that Schofield's result implies:

\begin{lem} For two dimension vectors ${\bf e}$ and ${\bf d}$, there exists an embedding $G_{\bf e}\subset G_{\bf d}$ if and only if ${\rm Ext}^1(G_{\bf e},G_{{\bf d}-{\bf e}})=0$.
\end{lem}

{\bf Proof:} Vanishing of ${\rm Ext}$ is equivalent to existence of an embedding $U\subset G_{\bf d}$, where $U$ is some representation of dimension vector ${\bf e}$. Thus $U$ is a degeneration of $G_{\bf e}$ and $G_{\bf d}/U$ is a degeneration of $G_{{\bf d}-{\bf e}}$. By the previous result, this yields a degeneration from $G_{{\bf d}-{\bf e}}* G_{\bf e}$ to $G_{\bf d}$, thus an isomorphism, proving that $G_{\bf e}$ embeds into $G_{\bf d}$.\\[1ex]
We can now formulate the main result of this section:

\begin{thm} The following are equivalent for a representation $M$ and a dimension vector ${\bf e}$ for the Dynkin quiver $Q$:
\begin{enumerate}
\item $M$ admits a subrepresentation of dimension vector ${\bf e}$,
\item $\dim{\rm Hom}(U,M)\geq\langle{\bf dim}\ U,{\bf e}\rangle$ for all indecomposable representations $U$ of $Q$,
\item $\dim{\rm Hom}(U,M)\geq\langle{\bf dim}\ U,{\bf e}\rangle$ for all representations $U$ of $Q$,
\item $\dim{\rm Hom}(U,M)\geq\langle{\bf dim}\ U,{\bf e}\rangle$ for all subrepresentations $U\subset G_{\bf e}$.
\end{enumerate}
\end{thm}

{\bf Proof:} The first statement implies the second as follows: if $N\subset M$ and ${\bf dim} N={\bf e}$, then
$$\langle{\bf dim}\, U,{\bf e}\rangle=\langle{\bf dim}\, U,{\bf dim}\, N\rangle=\dim{\rm Hom}(U,N)-\dim{\rm Ext}^1(U,N)\leq$$
$$\leq\dim{\rm Hom}(U,N)\leq\dim{\rm Hom}(U,M).$$

The implications from the second to the third and from the third to the fourth statement are trivial, and we will prove that the fourth statement implies the first by induction over ${\bf e}$, the case ${\bf e}=0$ being trivial.

The assumption of the fourth statement being fulfilled in particular for $U=G_{\bf e}$, we have
$$\dim{\rm Hom}(G_{\bf e},M)\geq\langle{\bf e},{\bf e}\rangle\geq 1$$
since $Q$ is Dynkin, thus there exists a non-zero map $f:G_{\bf e}\rightarrow M$. We define $K$, $I$ and $C$ as the kernel, image and cokernel of $f$, respectively, and denote ${\bf dim}\, K$ by ${\bf e}'$. Thus ${\bf e}'<{\bf e}$, and by induction, the theorem already holds for ${\bf e}'$.

Let $U$ be a subrepresentation of $G_{{\bf e}'}$. From ${\rm Ext}^1(G_{\bf e},G_{\bf e})=0$ we easily deduce ${\rm Ext}^1(K,I)=0$ since $Q$ is hereditary, and thus ${\rm Ext}^1(G_{{\bf e}'},G_{{\bf e}-{\bf e}'})=0$. By the above lemma, we find $G_{{\bf e}'}\subset G_{\bf e}$ and thus $U\subset G_{\bf e}$, which implies $\dim{\rm Hom}(U,M)\geq\langle{\bf dim}\, U,{\bf e}\rangle$ by assumption. We also easily derive ${\rm Ext}^1(U,I)=0$. Combination of these facts yields the estimate

$$\dim{\rm Hom}(U,C)\geq\dim{\rm Hom}(U,M)-\dim{\rm Hom}(U,I)=$$
$$=\dim{\rm Hom}(U,M)-\langle{\bf dim}\, U,{\bf dim}\, I\rangle\geq$$
$$\geq \langle {\bf dim}\, U,{\bf e}\rangle-\langle{\bf dim}\, U,{\bf dim}\, I\rangle=\langle{\bf dim}\, U,{\bf e}'\rangle.$$

This implies that the fourth assumption of the theorem is fulfilled for the dimension vector ${\bf e}'$ and the representation $C$; by the inductive assumption, there exists a subrepresentation $K'\subset C$ of dimension vector ${\bf e}'$. The inverse image of $K'$ under the projection $\pi$ in the exact sequence
$$0\rightarrow I\rightarrow M\stackrel{\pi}{\rightarrow}C\rightarrow 0$$
is thus a subrepresentation of $M$ of dimension vector ${\bf dim}\, I+{\bf dim}\, K'={\bf e}$, proving the theorem.\\[1ex]
{\bf Remark:} The second assumption of the theorem can easily be verified in practice, since it consists of finitely many inequalities which can be formulated entirely in terms of the Euler form of $Q$. The theorem can also be interpreted as an exact criterion for non-emptyness of the quiver Grassmannian ${\rm Gr}_{\bf e}(M)$.

\section{Irreducibility of quiver Grassmannians}\label{s5}

Again we assume $Q$ to be a Dynkin quiver. The following result is not immediately related to the previous sections, but arises naturally from the point of view of quiver Grassmannians:

\begin{thm} Given a dimension vector ${\bf e}$ and a representation $M$ as before, suppose the following inequalities holds:
\begin{enumerate}
\item $\dim{\rm Hom}(M,U)\leq\langle{\bf e},{\bf dim}\, U\rangle$ for all non-injective indecomposable $U$,
\item $\dim{\rm Hom}(U,M)\leq\langle{\bf dim}\, U,{\bf dim}\, M-{\bf e}\rangle$ for all non-projective indecomposables $U$.
\end{enumerate}

Then the quiver Grassmannian ${\rm Gr}_{\bf e}(M)$ is irreducible of dimension $$\dim{\rm Gr}_{\bf e}(M)=\langle{\bf e},{\bf dim}\, M-{\bf e}\rangle.$$
\end{thm}

{\bf Proof:} We divide the proof into three steps.
\begin{enumerate}
\item First, we prove that that the conclusion of the theorem holds if there exists an exact sequence
\begin{equation}\label{seq}0\rightarrow P\rightarrow M\rightarrow I\rightarrow 0\end{equation}
with $P$ projective, $I$ injective and ${\bf e}={\bf dim}\, P$. This is a slight generalization of \cite[Proposition 3.1]{CFR}, and we use the techniques developed in \cite[Section 2]{CFR}. First note that the quiver Grassmannian in question is non-empty since $P$ is a subrepresentation of $M$ of dimension vector ${\bf e}$. Moreover, $P$ also embeds into the representation of dimension vector ${\bf dim}\, M$ without selfextensions by \cite[Theorem 2.4]{Bo}. Thus, by \cite[Proposition 2.2]{CFR}, every irreducible component of ${\rm Gr}_{\bf e}(M)$ has at least dimension $\langle {\bf e},{\bf dim}\, M-{\bf e}\rangle$. If $C$ is such an irreducible component of ${\rm Gr}_{\bf e}(M)$, then $C$ is the closure of a set $\mathcal{S}_{[N]}$ for a representation $N$, where $\mathcal{S}_{[N]}\subset{\rm Gr}_{\bf e}(M)$ consists of all subrepresentations $U$ of $M$ which are isomorphic to $N$. The dimension of $\mathcal{S}_{[N]}$ equals $\dim{\rm Hom}(N,M)-\dim{\rm End}(N)$. Since $M$ degenerates to $P\oplus I$ and $P$ degenerates to $N$, we can estimate
$$\dim\mathcal{S}_{[N]}=\dim{\rm Hom}(N,M)-\dim{\rm End}(N)\leq$$
$$\leq \dim{\rm Hom}(N,P)-\dim{\rm End}(N)+\dim{\rm Hom}(N,I)\leq$$
$$\leq\dim{\rm Hom}(N,I)=\langle {\bf e},{\bf dim}\, M-{\bf e}\rangle.$$
If equality holds, then in particular $\dim{\rm Hom}(N,P)=\dim{\rm End}(N)$ and thus $N\simeq P$, again by \cite[Theorem 2.4]{Bo}. We thus conclude that $\mathcal{S}_{[P]}$ is the unique irreducible component of ${\rm Gr}_{\bf e}(M)$.
\item Second, we prove that existence of an exact sequence (\ref{seq}) is equivalent to the conditions
\begin{equation}\label{num}\dim{\rm Hom}(U,M)\leq\dim{\rm Hom}(U,P\oplus I),\end{equation}
$$\dim{\rm Hom}(M,U)\leq\dim{\rm Hom}(P\oplus I,U)$$ for all indecomposables $U$. Namely, by \cite[Theorem 3.3]{Bo} these numerical conditions are equivalent to existence of a degeneration of $M$ to $P\oplus I$. Again using \cite[Theorem 2.4]{Bo}, we conclude existence of the sequence (\ref{seq}).
\item Finally, we prove that the numerical conditions (\ref{num}) are implied by those of the theorem. If the conditions of the theorem hold, then, in particular, $\langle{\bf e},{\bf i}\rangle\geq 0$ for all vertices $i$ of $Q$, which implies that there exists a projective representation $P$ of dimension vector ${\bf e}$. Dually, we see that there exists an injective representation $I$ of dimension vector ${\bf dim}\, M-{\bf e}$. But then
$$\dim{\rm Hom}(U,M)\leq\langle {\bf dim}\, U,{\bf dim}\, M-{\bf e}\rangle=\langle {\bf dim}\, U,{\bf dim}\, I\rangle\leq$$
$$\leq\dim{\rm Hom}(U,I)\leq\dim{\rm Hom}(U,P\oplus I),$$
and similarly for the dual statement. This finishes the proof.
\end{enumerate}

\section{Appendix: A theorem of W. Crawley-Boevey and a linear algebra application}\label{app}

We review the main results of \cite{CB}. Let now $Q$ be an arbitrary finite quiver with set of vertices $Q_0$ and Euler form $\langle\_,\_\rangle$, let $M$ be a finite dimensional representation of $kQ$, and let ${\bf e}\in{\bf N}Q_0$ be a dimension vector for $Q$.\\[1ex]
We view ${\bf e}$ as a functional on the Grothendieck group of $kQ$ via ${\bf e}(M)=\langle\dimv M,{\bf e}\rangle$. The representation $M$ is called ${\bf e}$-semistable if ${\bf e}(M)=0$ and ${\bf e}(N)\geq 0$ for all subrepresentations $N\subset M$.\\[1ex]
We denote by ${\rm hom}(M,{\bf e})$ the generic (that is, minimal possible) dimension of the space of homomorphisms from $M$ to a representation of dimension vector ${\bf e}$, and by ${\bf f}_{M,{\bf e}}\in{\bf N}Q_0$ the generic (that is, maximal possible) rank of a map from $M$ to a representation of dimension vector ${\bf e}$. Similarly to the above, we denote by ${\rm Gr}^{\bf e}(M)$ the Grassmannian of subrepresentations of $M$ of {\it codimension} ${\bf e}$. Then the following holds:
\begin{thm}\label{cb}\cite[Theorem 1,Corollary 1,Corollary 2, Proof of Corollary 1]{CB}~
\begin{enumerate}
\item We have $${\rm hom}(M,{\bf e})=\langle {\bf f}_{M,{\bf e}},{\bf e}\rangle+\dim U$$
for some nonempty open subset $U$ of ${\rm Gr}^{{\bf f}_{M,{\bf e}}}(M)$.
\item We have
$$\lim_{r\rightarrow\infty}\frac{1}{r}{\rm hom}(M,r\cdot{\bf e})=\max\{{\bf e}(M/N)\, :\, N\subset M\}.$$
More precisely, the sequence $(\frac{1}{r}{\rm hom}(M,r\cdot{\bf e}))_r$ converges to its infimum, which is the right hand side.
\item If $M$ is ${\bf e}$-semistable, then ${\rm hom}(M,r\cdot{\bf e})=0$ for some $r\geq 1$.
\end{enumerate}
\end{thm}

We will use these results (and some of the methods of their proofs) to deduce:

\begin{thm} If ${\bf e}(N)\geq 0$ for all subrepresentations $N\subset M$ , then $${\rm hom}(M,r\cdot{\bf e})=r\cdot{\bf e}(M)$$ {\it for all large enough} $r$.
\end{thm}

{\bf Proof:} First we reduce to the case where ${\bf e}(N)>0$ for all non-zero proper subrepresentations $N\subset M$. Namely, suppose that ${\bf e}(M')=0$ for a non-zero subrepresentation $M'\subset M$, and consider the sequence
$$0\rightarrow M'\rightarrow M\stackrel{\pi}{\rightarrow}\bar{M}\rightarrow 0.$$
If $N\subset M'$, then ${\bf e}(N)\geq 0$, thus $M'$ is ${\bf e}$-semistable. If $N\subset \bar{M}$, then $\pi^{-1}N\subset M$, thus
$${\bf e}(N)={\bf e}(\pi^{-1}N)-{\bf e}(M')={\bf e}(\pi^{-1}N)\geq 0,$$
thus $\bar{M}$ still satisfies the assumption of the theorem. If the theorem holds for $M'$ and $\bar{M}$, then, for all large enough $r$,  we have
$${\rm hom}(M,r\cdot{\bf e})\leq{\rm hom}(M',r\cdot{\bf e})+{\rm hom}(\bar{M},r\cdot{\bf e})=$$
$$={\rm hom}(\bar{M},r\cdot{\bf e})=r\cdot{\bf e}(\bar{M})=r\cdot{\bf e}(M),$$
and thus equality by the second part of Theorem \ref{cb}. Now an induction over the dimension reduces the theorem to the case where ${\bf e}(N)>0$ for all non-zero proper $N\subset M$.\\[1ex]
In this case, the second part of the theorem reads
$$\lim_{r\rightarrow\infty}\frac{1}{r}{\rm hom}(M,r\cdot{\bf e})={\bf e}(M),$$
and, by the first part,
$${\rm hom}(M,r\cdot{\bf e})=r\cdot{\bf e}({\bf f}_{M,r\cdot{\bf e}})+\dim U$$
for some non-empty open $U\subset{\rm Gr}^{{\bf f}_{M,r\cdot{\bf e}}}(M)$. Since ${\bf f}_{M,r\cdot{\bf e}}$ can assume only finitely many values, at least one value ${\bf d}$ is assumed infinitely many times. Let us fix an arbitrary ${\bf d}$ with this property; thus there exists a subsequence $r_1<r_2<\ldots$ of the integers such that
$${\rm hom}(X,r_i\cdot{\bf e})=r_i\cdot{\bf e}({\bf d})+\dim U_i$$
for all $i\geq 1$ and varieties $U_i$ of bounded dimensions. Thus
$$\frac{1}{r_i}{\rm hom}(M,r_i\cdot{\bf e})={\bf e}({\bf d})+\frac{1}{r_i}\dim U_i$$
for all $i\geq 1$. Taking limits of both sides, we see that
$${\bf e}(M)=\lim_{i\rightarrow\infty}\frac{1}{r_i}{\rm hom}(M,r_i\cdot{\bf e})={\bf e}({\bf d})+\lim_{i\rightarrow\infty}\frac{1}{r_i}\dim U_i={\bf e}({\bf d}).$$
Since ${\bf e}({\bf d})={\bf e}(M/N)$ for some subrepresentation $N\subset M$, we have
$${\bf e}(M)={\bf e}(M/N)={\bf e}(M)-{\bf e}(N),$$
which, by assumption, is only possible if $N=0$ or $N=M$, thus ${\bf d}=\dimv M$ or ${\bf d}=0$. In either case, ${\rm Gr}^{\bf d}(M)$ reduces to a point, thus $U_i$ consists of a single point, and thus ${\rm hom}(M,r\cdot{\bf e})$ equals $r\cdot{\bf e}(M)$ for all the (infinitely many) values $r$ such that ${\bf f}_{M,r\cdot{\bf e}}={\bf d}$. This conclusion applying to all ${\bf d}$ which occur as ${\bf f}_{M,r\cdot{\bf e}}$ infinitely often, we see that failure of ${\rm hom}(M,r\cdot{\bf e})=r\cdot{\bf e}(M)$ can only happen for finitely many $r$, proving the theorem.\\[1ex]
From the previous theorem, we will derive the following linear algebra result:

\begin{lem}\label{lalapp} Let $Z\subset{\rm Hom}(V,W)$ be a subspace of linear maps between finite dimensional $k$-vector spaces $V$ and $W$. Assume that $\dim Z(U)\geq\dim U$ for all $U\subset V$. Then, for all large enough $r$, there exists a matrix $F\in M_{r\times r}(Z)$ defining an injective map from $V^r$ to $W^r$.
\end{lem}

{\bf Proof:} Choose a basis $f_1,\ldots,f_k$ of $Z$ and define a representation $M$ of the $k$-arrow Kronecker quiver $Q=K_k$ (with $k$ arrows from vertex $1$ to vertex $2$) by $f_1,\ldots,f_k:V\rightarrow W$. We consider the dimension vector ${\bf e}=(k-1,1)$.  For a representation $N$ of $K_k$, we then have ${\bf e}(N)=\dim Y_2-\dim Y_1$. We claim that ${\bf e}(N)\geq 0$ for all $N\subset M$.
Namely, if $N$ is a subrepresentation of $M$, then $N_1$ is a subspace of $M_1=V$ and $N_2$ is a subspace of $M_2=W$ containing $Z(N_1)$, thus, by assumption,
$$\dim N_2\geq \dim Z(N_1)\geq \dim N_1,$$
which implies
$${\bf e}(N)=\dim N_2-\dim N_1\geq 0.$$
By the above theorem, we thus have
$${\rm hom}(M,r\cdot{\bf e})=r\cdot(\dim W-\dim V)$$
for all large enough $r\geq 1$. A general representation $U$ of dimension vector $r\cdot{\bf e}$ can be written as the kernel of a general map between injectives
$$0\rightarrow U\rightarrow I_1^r\rightarrow I_2^r\rightarrow 0.$$
Applying ${\rm Hom}(M,\_)^*$, this yields an exact sequence
$$\underbrace{{\rm Hom}(M,I_1^r)^*}_{\simeq V^r}\rightarrow\underbrace{{\rm Hom}(M,I_2^r)^*}_{\simeq W^r}\rightarrow{\rm Hom}(M,U)^*\rightarrow 0.$$
Additionally choosing $U$ such that $\dim{\rm Hom}(M,U)=r\cdot(\dim W-\dim V)$, we see that the map $F:V^r\rightarrow W^r$ in this sequence, which is given by a matrix with entries being linear combinations of the $f_i$, or, in other terms, $F\in M_{r\times r}(Z)$, is injective.

\section*{Acknowledgements} The second named author would like to thank K. Bongartz for many inspiring discussions on Theorem \ref{thm1}, W. Crawley-Boevey for useful remarks on Section \ref{section_ex}, and C. M. Ringel for  an example in Section \ref{section_ex}.

\end{document}